\documentclass[twoside]{amsart}
\usepackage{amsmath,amsthm,amssymb}
\input{xy}
\xyoption{all}
\def\bar{\overline}
\def\today{Jan 16, 2004}

\theoremstyle{plain}
\newtheorem{theorem}{Theorem}[section]
\newtheorem{proposition}[theorem]{Proposition}
\newtheorem{lemma}[theorem]{Lemma}

\newtheorem{corollary}[theorem]{Corollary}

\theoremstyle{definition}

\newtheorem{remark}[theorem]{Remark}
\theoremstyle{remark}
\newtheorem{acknowledgments}{Acknowledgments}

\newcommand{\bA}{{\bf A}}
\newcommand{\bB}{{\bf B}}

\newcommand{\bM}{{\bf M}}
\newcommand{\bN}{{\bf N}}
\newcommand{\bP}{{\bf P}}

\newcommand{\A}{{\mathbb A}}
\newcommand{\F}{{\mathbb F}}

\newcommand{\Q}{{\mathbb Q}}
\newcommand{\R}{{\mathbb R}}
\newcommand{\Z}{{\mathbb Z}}

\newcommand{\cC}{ {\mathcal C} }
\newcommand{\cH}{ {\mathcal H} }
\newcommand{\cI}{ {\mathcal I} }
\newcommand{\cJ}{ {\mathcal J} }

\newcommand{\cO}{ {\mathcal O} }
\newcommand{\cP}{ {\mathcal P} }

\newcommand{\hP}{\hat{P}}

\newcommand{\lcm}{ {\rm lcm } }
\newcommand{\ord}{ {\rm ord} }
\newcommand{\Gal}{{\rm Gal}}

\newcommand{\NP}{ {\rm  NP} }     
\newcommand{\HP}{ {\rm  HP} }     
\newcommand{\Tr}{ {\rm Tr} }
\newcommand{\N}{ {\rm N} }

\newcommand{\Zeta}{{\rm Zeta}}
\renewcommand{\det}{{\rm det}}

\begin{document}

\title[$L$-functions of Exponential sums over one-dimensional affinoids]
{ $L$-functions of exponential sums over one-dimensional
affinoids: Newton over Hodge}
\author{Hui June Zhu}
\address{
Department of mathematics and statistics,
McMaster University,
Hamilton, ON L8S 4K1,
CANADA
}
\email{zhu@cal.berkeley.edu}

\date{\today}
\subjclass[2000]{11,14}

\begin{abstract}
This paper proves a sharp lower bound for Newton polygons of
$L$-functions of exponential sums of one-variable rational
functions. Let $p$ be a prime and let $\bar{\F}_p$ be the
algebraic closure of the finite field of $p$ elements. Let
$\bar{f}(x)$ be any one-variable rational function over
$\bar{\F}_p$ with $\ell$ poles of orders $d_1,\ldots,d_\ell$.
Suppose $p$ is coprime to $d_1\cdots d_\ell$. We prove that there
exists a tight lower bound which we call Hodge polygon, depending
only on the $d_j$'s, to the Newton polygon of $L$-function of
exponential sums of $\bar{f}(x)$. Moreover, we show that for any
$\bar{f}(x)$ these two polygons coincide if and only if $p\equiv
1\bmod d_j$ for every $1\leq j\leq \ell$. As a corollary, we
obtain a tight lower bound for the $p$-adic Newton polygon of
zeta-function of an Artin-Schreier curve given by affine equations
$y^p-y=\bar{f}(x)$.
\end{abstract}

\maketitle

\section{Introduction}\label{S:1}

Let $\A$ be the space of rational functions in one variable $x$
with $\ell$ distinct poles
(say at $P_1, P_2,\ldots, P_\ell$)
of orders $d_1, \ldots, d_\ell\geq 1$ on the projective line.
For any field $K$, we denote by $\A(K)$
the set of all rational
functions of the form
$
\sum_{j=1}^{\ell}\sum_{i=1}^{d_j} a_{j,i} (x-P_j)^{-i},
$
where coefficients $a_{j,i}\in K$, poles $P_j\in K\cup \{\infty\}$
and $\prod_{j=1}^{\ell}a_{j,d_j}\neq 0$ (we set
$(x-\infty)^{-i}=x^i$ for the point at $\infty$). Naturally one
may consider $\A$ as a quasi-affine space parameterized by
coefficients $a_{j,i}$ for all $i\geq 1,1\leq j\leq \ell$ and
poles $P_j$ for $1\leq j\leq \ell$. Let the {\em Hodge polygon} of
$\A$, denoted by $\HP(\A)$, be the end-to-end join of line
segments of horizontal length $1$ with slopes listed below:
\begin{equation}\label{E:HP}
\overbrace{0,\ldots,0}^{\ell-1};
\overbrace{1,\ldots,1}^{\ell-1};
\overbrace{\frac{1}{d_1}, \cdots, \frac{d_1-1}{d_1}}^{d_1-1};
\overbrace{\frac{1}{d_2}, \cdots, \frac{d_2-1}{d_2}}^{d_2-1};
\ldots\ldots;
\overbrace{\frac{1}{d_\ell}, \cdots, \frac{d_\ell-1}{d_\ell}}^{d_\ell-1}.
\end{equation}
They are joined in a nondecreasing order from left to right
starting from the origin on $\R^2$.
Let $d:=\sum_{j=1}^{\ell}d_j+\ell-2$. So $\HP(\A)$ is a
lower convex hull in $\R^2$ with endpoints $(0,0)$ and $(d,d/2)$.

For any prime $p$ let
$E(x)=\exp(\sum_{i=0}^{\infty}\frac{x^{p^i}}{p^i})$ be the
$p$-adic Artin-Hasse exponential function. Let $\gamma$ be a $p$-adic
root of $\log(E(x))$ with $\ord_p \gamma=\frac{1}{p-1}$.
Then $E(\gamma)$ is a primitive $p$-th root of unity.
We fix it and denote it by $\zeta_p$.

In this paper we let $p$ be a prime
coprime to $\prod_{i=1}^{\ell}d_i$ and let $a$ be
a positive integer. Let $q=p^a$.
Let
\begin{eqnarray}\label{E:fbar}
\bar{f}(x) &: =& \sum_{j=1}^{\ell}\sum_{i=1}^{d_j}
\bar{a}_{j,i} (x-\bar{P}_j)^{-i},
\end{eqnarray}
where
$\bar{a}_{j,i}\in \F_q, \bar{P}_j\in\F_q\cup\{\infty\}$
for every $i,j$.
Let $\bar{g}(x):=\prod_{\bar{P}_j\neq \infty}(x-\bar{P}_j)
\in\F_q[x]$.
For any positive integer $k$, the $k$-th exponential sum of
$\bar{f}(x)\in\F_q(x)$ is
$S_k(\bar{f}):=
\sum
\zeta_p^{\Tr_{\F_{q^k}/\F_p}(\bar{f}(x))}
$
where the sum ranges over all $x$ in $\F_{q^k}$ such that
$\bar{g}(x)\neq 0$.
The $L$-function of
the exponential sum of $\bar{f}$ is defined by
\begin{eqnarray}\label{E:Lfunction}
L(\bar{f};T):
=\exp(\sum_{k=1}^{\infty}S_k(\bar{f})\frac{T^k}{k}).
\end{eqnarray}
It is well known that the $L$-function is a polynomial in
$\Z[\zeta_p][T]$ of
degree $d$ (e.g., by
combining the Weil Conjecture for curves with the argument in
\cite{Zhu:1} between Remark 1.2 and Corollary 1.3). One may write
\begin{eqnarray}\label{E:20}
L(\bar{f};T)&=&1+b_1T+b_2T^2+\ldots+b_{d}T^{d}\in\Z[\zeta_p][T].
\end{eqnarray}
Define the {\em Newton polygon} of the $L$-function of
$\bar{f}$ over $\F_q$ as the lower convex hull in $\R^2$
of the points $(n,\ord_q(b_n))$ with $0\leq n\leq d$, where we set
$b_0=1$ and $\ord_q(\cdot):=\ord_p(\cdot)/a$.
We denote it by $\NP(\bar{f}; \F_q)$. One notes
immediately that the Newton polygon $\NP(\bar{f};\F_q)$ and the
Hodge polygon $\HP(\A)$ have the same endpoints $(0,0)$ and
$(d,d/2)$.
Let $\lcm(d_j)$
denote the least common multiple
of $d_j$'s for all $1\leq j\leq \ell$.
The main result of the present paper is the following.

\begin{theorem}\label{T:main1}
Let notation be as above.
For any rational function $\bar{f}\in\A(\F_q)$,
the  Newton polygon $\NP(\bar{f};\F_q)$
lies over the Hodge polygon $\HP(\A)$, and their endpoints meet.
Moreover, for any $\bar{f}\in \A(\F_q)$
one has $\NP(\bar{f};\F_q)=\HP(\A)$
if  and only if $p\equiv 1\bmod (\lcm\; d_j)$.
\end{theorem}

\begin{remark}\label{R:main1}
The first part (i.e., Newton over Hodge) of Theorem \ref{T:main1} was
a conjecture of Adolphson-Sperber and Bjorn Poonen, described to the author
independently in 2001. The case $\ell=1$
is known (see \cite{Sperber} or \cite{Wan:1}). The case $\ell=2$
and $\bar{f}(x)$ has only poles at $\infty$ and $0$ (i.e.,
$\bar{f}(x)$ is a one variable Laurent polynomial) was obtained first
by Robba (see \cite[Theorems 7.2 and 7.5]{Robba:1}).  Theorem
\ref{T:main1} is an analog of Katz-type conjectures (see \cite[Theorem
2.3.1]{Katz} and \cite{Mazur}).
\end{remark}

Below we shall discuss some
applications of our result in algebraic geometry.
A question that remains open is
whether there is a curve in every Newton polygon stratus in
the moduli space of curves
over $\bar{\F}_p$ (for every $p$).
Recently \cite{Geer-Vlugt:92}
and \cite{Geer-Vlugt:95} gave an affirmative answer
to this question for $p=2$ by constructing supersingular curves over
$\bar{\F}_2$ via a fibre product of Artin-Schreier curves.
It is essential to understand the shape of
Newton polygons of Artin-Schreier
curves, and in particular, to find a sharp lower bound for them.
The Newton polygon of the Artin-Schreier curve
$C_{\bar{f}}: y^p-y=\bar{f}(x)$ over $\F_q$
(note that $C_{\bar{f}}$ has genus $d(p-1)/2$)
is the normalized $p$-adic Newton polygon of the numerator of
the Zeta function $\Zeta(C_{\bar{f}};T)$ of $C_{\bar{f}}$
(here `normalized' means taking $\ord_p(\cdot)/a$ as the
valuation).
Denote this Newton polygon by $\NP(C_{\bar{f}};\F_q)$. Then we have
the following corollary.

\begin{corollary}\label{C:main1}
Let notation be as in Theorem \ref{T:main1}.
For any $\bar{f}\in\A(\F_q)$ and Artin-Schreier curve
$C_{\bar{f}}:y^p-y=\bar{f}(x)$, the Newton polygon
$\NP(C_{\bar{f}};\F_q)$ shrunk by a factor of $1/(p-1)$
(vertically and horizontally)
is equal to $\NP(\bar{f};\F_q)$ and it lies over the Hodge
polygon $\HP(\A)$.
Moreover, for any $\bar{f}\in \A(\F_q)$
the equality holds
if and only if $p\equiv 1\bmod (\lcm\; d_j)$.
\end{corollary}
\begin{proof}
We shall first give an elementary proof of the following relation between
the Zeta function of $C_{\bar{f}}$ and the $L$-function of $\bar{f}$:
\begin{equation}\label{E:zeta}
\Zeta(C_{\bar{f}};T) = \frac{\N_{\Q(\zeta_p)/\Q}(L(\bar{f};T))}{(1-T)(1-qT)}
\end{equation}
with the norm $\N_{\Q(\zeta_p)/\Q}(\cdot)$
being interpreted as the product of conjugates of the $L$-function
$L(\bar{f};T)$
in $\Q(\zeta_p)$ over $\Q$, the automorphism acting trivially
on the variable $T$.
(One may also see, for example,
\cite[Section VI, (93)]{Bombieri} for some relevant discussion.)
Recall that for any integer $n$ one has
$\sum_{a\in\F_p}\zeta_p^{an}=p$ or $0$
depending on whether $n$ is $0$ or not, respectively.
For any $k\geq 1$ let
$C_{\bar{f}}'$ be the curve $C_{\bar{f}}$
less the $\ell$ ramification points over
$\bar{P}_1=\infty,\bar{P}_2, \ldots,\bar{P}_\ell$.
Write
$\F_{q^k}^+:=\F_{q^k}-\{\bar{P}_2,\ldots,\bar{P}_\ell\}$.
Then
$$\# C_{\bar{f}}'(\F_{q^k}) = \sum_{a\in\F_p}\sum_{x\in\F_{q^k}^+}
\zeta_p^{a\Tr(\bar{f}(x))}
$$
where $\Tr(\cdot) = \Tr_{\F_{q^k}/\F_p}(\cdot)$.
It follows that
\begin{eqnarray*}
\Zeta(C_{\bar{f}};\F_q)
&=&\exp\left(\sum_{k=1}^{\infty}
       \left(\ell+\#C_{\bar{f}}'(\F_{q^k})\right)\frac{T^k}{k}\right)\\
&=&\exp\left(\sum_{k=1}^{\infty}(1+q^k+\sum_{a\in\F_p^*}
   \sum_{x\in\F_{q^k}^+}\zeta_p^{a\Tr(\bar{f}(x))})\frac{T^k}{k}\right)\\
&=&\frac{
   {\prod_{a\in\F_p^*}\left(\exp(\sum_{k=1}^{\infty}
   (\sum_{x\in\F_{q^k}^+}{\zeta_p^{a\Tr(\bar{f}(x))}}))\frac{T^k}{k}\right)}}
     {(1-T)(1-q T)}\\
&=&\frac{\N_{\Q(\zeta_p)/\Q}(L(\bar{f};T))}{(1-T)(1-q T)}.
\end{eqnarray*}
This proves (\ref{E:zeta}). Our first assertion of the corollary
follows from (\ref{E:zeta}) and the rest then follows from the
first assertion and Theorem \ref{T:main1}.
\end{proof}

\begin{remark}\label{R:main2}
We remark that $\NP(\bar{f};\F_q)$ and $\HP(\A)$ always coincide at
the slope-0 segments (of horizontal length $\ell-1$).  Indeed, in the
spirit of Corollary \ref{C:main1}, it suffices to show
that the Artin-Schreier curve $C_{\bar{f}}: y^p-y=\bar{f}(x)$ has
$p$-rank equal to $(\ell-1)(p-1)$, which follows from {\em
Deuring-Shafarevic formula} (see, for instance, \cite[Corollary
1.8]{Crew}). By symmetry, their slope-$1$ segments also coincide.
\end{remark}

Now we give an amusing example:
By the above, the curve
$$
C/\bar{\F}_p: y^p-y = a_{1,2}x^2 + a_{1,1}x +
\sum_{j=2}^{\ell}(\frac{a_{j,1}}{(x-\bar{P}_j)^2}+\frac{a_{j,2}}{(x-\bar{P}_j)})
$$
for odd prime $p$ and nonzero $a_{j,i}$'s,
has its Newton polygon
slope $0$ (resp. $1$) of length $(\ell-1)(p-1)$,
and slope $1/2$ of length $\ell(p-1)$.

Finally we comment on
our conventions for the proof of this theorem:
We first note that for
$d_1=\ell=1$ the $L$-function of $f$ is equal to $1$ and we shall
exclude this case for the rest of the paper for simplicity;
Following conventions often used in algebraic geometry,
we set $P_1=\infty$. Recall that one
can always move one pole of an Artin-Schreier curve
to $\infty$ by an automorphism of the projective line it covers
without altering its zeta function.
If $\ell>1$ then we set $P_2=0$ for the rest of the paper.
This is not a restriction for our
purpose since it is observed by definition that
$L(\bar{f}(x+c);T)=L(\bar{f}(x);T)$ for any $c\in \F_q$ so one can
always shift $\bar{f}(x)$ so that one of its poles lies at $0$.
(Note that we then have $\bar{g}(0)=0$.)

This paper is organized as follows. We develop core theory of
exponential sums over an affinoid in Section 2. There we determine
the size of residue disks and derive an effective trace formula of
Dwork-Monsky-Reich. In other words, Section \ref{S:2} contains
fundamentals for the rest of the paper. In Section \ref{S:3}, we
present a practical algorithm to estimate the $p$-adic valuation
of our Frobenius matrix. Finally Section \ref{S:4} is devoted to
the proof of Theorem \ref{T:main1}. In a sequel paper
\cite{Li-Zhu:1}
we shall study
the asymptotic
(as $p$ varies) generic Newton polygons for
$L$-functions of exponential sums of
one-variable rational function.

\section{Exponential sums over one dimensional affinoids}\label{S:2}

We generalize Robba's work in \cite{Robba:1} from one dimensional
annuli to one dimensional affinoids. This differs from Dwork's
original approach (see \cite{Dwork1,Dwork2}), but is somewhat akin
to \cite{Berthelot}. Our exposition is (of course) after
\cite{Monsky-Washnitzer:I, Robba:1,Robba:2,Monsky:II, Monskey:III,
Reich}. For fundamental material considering non-Archimedean
geometry see \cite{BGR} or \cite{Schikhop}.

\subsection{Preliminaries}\label{S:2.0}

Let $\Q_q$ be the degree
$a$ unramified extension of $\Q_p$ and let $\Z_q$ be its ring of
integers. Let $\bar\Q_p$ be the algebraic closure of $\Q_p$, and
let $\bar\Z_p$ its ring of integers.
Let $\Omega$ be the $p$-adic completion of $\bar\Q_p$.
Let $\Omega_1=\Q_p(\zeta_p)$ and $\Omega_a$ the unique unramified
extension of $\Omega_1$ of degree $a$ in $\Omega$. Let $\cO_1$ and
$\cO_a$ be the rings of integers in $\Omega_1$ and $\Omega_a$,
respectively.
Note that $\cO_1=\Z_p[\zeta_p]=\Z_p[\gamma]$.
Fix roots $\gamma^{1/d_j}$ in $\cO_a$
for the rest of the paper, and
let $\Omega'_1:=\Omega_1(\gamma^{1/d_1}, \cdots, \gamma^{1/d_\ell})$.
Let $\Omega'_a:=\Omega'_1\Omega_a$. Let $\cO'_a$ and $\cO'_1$ be
the rings of integers of $\Omega'_a$ and $\Omega'_1$, respectively.
 Let $|\cdot|_p$ be the $p$-adic valuation on
$\Omega_a$ such that $|p|_p=p^{-1}$.
The diagrams below represent these field extensions and
the associated extensions of their rings of integers.
$$
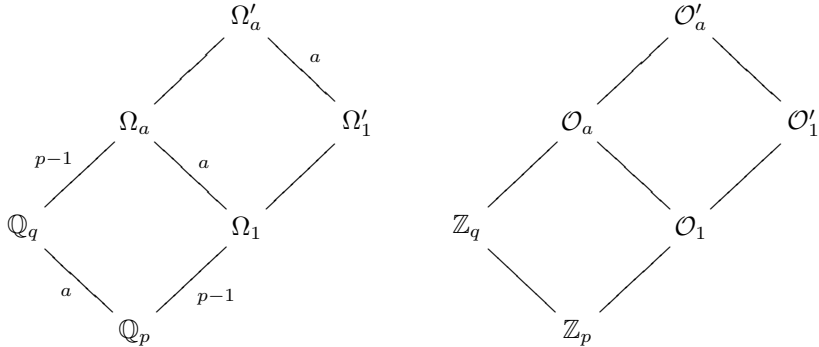
\xymatrix{
 &&&{\Omega'_a}\ar @{-}[dl]\ar @{-}[dr]^{a} &&&
   &{\cO'_a}\ar @{-}[dl]\ar @{-}[dr]&\\
  &&{\Omega_a}\ar @{-}[dl]_{p-1}\ar @{-}[dr]^{a}
   && {\Omega'_1}\ar @{-}[dl]
   &&{\cO_a}\ar @{-}[dl] \ar @{-}[dr]
   && {\cO'_1}\ar @{-}[dl]\\
&{\Q_q} \ar @{-}[dr]_{a}
  &&{\Omega_1} \ar @{-}[dl]^{p-1}
  &&{\Z_q} \ar @{-}[dr]
  &&{\cO_1} \ar @{-}[dl]
  &
  \\
&&{\Q_p}
  &
  &
  && {\Z_p}
  &&
}
$$

By taking Teichm\"uller lifts of coefficients and poles of
$\bar{f}\in \F_q[x]$, we get $\hat{f}(x)\in\Z_{q}[x]$ with
$\hat{f}(x) = \sum_{i=1}^{d_1}a_{1,i}x^{i} +
\sum_{j=2}^{\ell}\sum_{i=1}^{d_j}a_{j,i} (x-\hP_j)^{-i}$. Since
$\hP_j$ is a Teichm\"uller lift one has $\hP_j^\tau = \hP_j^p$.
Similarly, let $\hat{g}(x)=\prod_{j=2}^{\ell}(x-\hP_j)\in\Z_q[x]$ be
the corresponding Teichm\"uller lift of $\bar{g}(x)\in\F_q[x]$.
Note that $\hP_j^q=\hP_j$ and $\hP_j^q\equiv \bar{P}_j\bmod \cP$.

We mainly work on $p$-adic spaces over $\Omega_a$ (or
$\Omega'_a$). Let $|\Omega_a|_p$ denote the $p$-adic value group
of $\Omega_a$. Let $\bP^1$ be the rigid projective line over
$\Omega_a$. For any $\hP\in\Omega_a$ and $r\in |\Omega_a|_p$ let
$\bB[\hP,r]$ and $\bB(\hP,r)$ denote the closed disk and (wide) open
disk of radius $r$ about $\hP$ on $\bP^1$, that is
$\bB[\hP,r]:=\{X\in\Omega_a| |X-\hP|_p\leq r\}$ and
$\bB(\hP,r):=\{X\in\Omega_a| |X-\hP|_p<r\}$.

Let $r\in |\Omega_a|_p$ and $0<r<1$. For any positive $p$-power
$s$ let $\bA_{r,s} :=
\bB[0,1/r]-\bigcup_{j=2}^{\ell}\bB(\hP_j^s,r) =
\bP^1-\bigcup_{j=1}^{\ell}\bB(\hP_j^s,r)$. So $\bA_{r,s}= \{X\in
\Omega_a| |X|_p\leq 1/r; |X-\hP_j^s|_p\geq r \mbox{ for } 2\leq
j\leq \ell,\}$ and it is an affinoid over $\Omega_a$. Let
$\cH(\cdot)$ be the ring of rigid analytic functions over
$\Omega_a$ of a given affinoid. Hence, $\cH(\bA_{r,s})$ is a
$p$-adic Banach space with the natural $p$-adic supremum norm. For
ease of notation, we shall abbreviate $\bA_r$ for $\bA_{r,1}$ in
this paper.

\subsection{The $p$-adic Mittag-Leffler decomposition}\label{S:2.1}

Let $\cH_1(\bA_r)$ be the subset of $\cH(\bA_r)$ consisting of all
rigid analytic functions on $\bB[0,1/r]$. For $2\leq j\leq \ell$
let $\cH_j(\bA_r)$ be the subset of $\cH(\bA_r)$ consisting of all
rigid analytic functions on $\bP^1-\bB(\hP_j,r)$ that are
holomorphic at $\infty$ and vanish at $\infty$. For any rigid
analytic function $\xi$ defined on a subset $B$ of $\bP^1$, let
$||\xi||_B:=\sup_{x\in B}|\xi(x)|_p$ (i.e., supremum norm). This
defines a norm on $\cH(\bA_r)$ and $\cH_j(\bA_r)$, which are $p$-adic
Banach spaces under the supremum norm $||\xi ||_{\bA_r}$.

\begin{lemma}[$p$-adic Mittag-Leffler]\label{L:KML}
Let $r\in |\Omega_a|_p$ and $0<r<1$.
Then the $\bB(\hP_j,r)$'s with $1\leq j\leq \ell$ are mutually
disjoint. There is a canonical decomposition of $p$-adic Banach
spaces $\cH(\bA_r)\cong \bigoplus_{j=1}^{\ell}\cH_j(\bA_r )$ in the
sense that for any $\xi\in\cH(\bA_r)$ there is a unique
$\xi_{\hP_j}\in \cH_j(\bA_r)$ such that every  $\xi-\xi_{\hP_j}$ is
analytically expandable to $\bB(\hP_j,r)$. Every $\xi$ can be
uniquely represented as a sum $\xi=\sum_{j=1}^{\ell}\xi_{\hP_j}$
such that
\begin{eqnarray}\label{E:cH2}
||\xi||_{\bA_r} & = & \max_{1\leq j\leq
\ell}(||\xi_{\hP_j}||_{\bP^1-\bB(\hP_j,r)}).
\end{eqnarray}
\end{lemma}
\begin{proof}
We first show the disjointness.
Let $j\geq 3$.
Since the $\hP_j$'s are Teichm\"uller lifts in $\Z_q$ with
$\hP_j^q=\hP_j$ one has $|\hP_j|_p=1$.
For any $3\leq i<j\leq \ell$,
one first observes easily
that $|\hP_i-\hP_j|_p\leq \max(|\hP_i|_p,|\hP_j|_p)=1$.
The hypothesis that $\bar{P}_i\neq \bar{P}_j$ in
the residue field of $\Z_q$ implies that
$|\hP_i-\hP_j|_p\nless 1$ and hence one has
$|\hP_i-\hP_j|_p=1$.

Let $j\geq 2$. Pick any $\hP\in \bB(\hP_j,r)$. If $j=2$ then
$|\hP|_p<r<1<1/r$ so $\hP\in \bB[1,1/r]$; If $j\geq 3$ then
$|\hP-\hP_j|_p<r<1$ and $|\hP_j|_p=1$ imply that $|\hP|_p=1<1/r$ so $\hP\in
\bB[0,1/r]$. This proves that $\bB(\infty,r)\cap
\bB(\hP_j,r)=\emptyset$ for all $j\geq 2$.

Now let $j,j'\geq 2$ and $j'\neq j$. For any $\hP\in \bB(\hP_j,r)$,
one has $|\hP-\hP_j|_p<r<1$ and $|\hP_j-\hP_{j'}|_p=1$ by the previous
paragraph and so $|\hP-\hP_{j'}|_p=1>r$. This shows $\hP\not\in
\bB(\hP_{j'},r)$. This proves the disjointness.

The proof for the rest of the lemma follows directly from
\cite[Theorem 4.7]{Robba:0} (see also \cite{Krasner}).
\end{proof}

Every element in the $\Omega_a$-space
$\cH(\bA_r)$ can be uniquely represented as
$
\sum_{i\geq 0}c_{1,i}X^i + \sum_{j=2}^{\ell}\sum_{i\geq
1}c_{j,i}(X-\hP_j)^{-i}
$
where $c_{j,i}\in\Omega_a$ and $\forall j\geq 1,
 \lim_{i\rightarrow\infty}\frac{|c_{j,i}|_p}{r^i}=0$.
For simplicity, we write $X_1:=X$ and $X_j:=(X-\hP_j)^{-1}$
for $1\leq j\leq \ell$.
Then the $\Omega_a$-space $\cH(\bA_r)$
has a natural monomial basis
$\vec{b}_{\rm unw}:=\{1,X_1^i,X_2^i,\ldots,X_\ell^i\}_{i\geq 1}.$
In Theorem \ref{T:reduce1}
we shall use a {\em weighted} basis $\vec{b}_{\rm w}$
(note that neither $\vec{b}_{\rm unw}$ nor $\vec{b}_{\rm w}$
is an orthonormal basis):

\begin{lemma}\label{L:basis}
Let $r\in |\Omega_a|_p$ and $0<r<1$.
For each $1\leq j\leq \ell$, let
$Z_j:=\gamma^{\frac{1}{d_j}}X_j$.
Then $\vec{b}_{\rm w}:=\{1,Z_1^i,\ldots,Z_\ell^i\}_{i\geq 1}$
forms a basis
of $\cH(\bA_r)$ over $\Omega'_a$.
\end{lemma}
\begin{proof}
Obvious by Lemma \ref{L:KML} and remarks
preceding the lemma.
\end{proof}

\subsection{The $U_p$ operator}\label{S:2.2}

For any $s\in p^{\Z_{\geq 0}}$ and for any
$\xi(X)\in\cH(\bA_{r,s})$, let $U_p$ be the map defined by $(U_p
\xi)(X) := 1/p \cdot\sum_{Z^p=X}\xi(Z)$ from $\cH(\bA_{r,s})$ to
$\cH(\bA_{r^p,sp})$. Similarly let $(U_q  \xi)(X) := 1/q
\cdot\sum_{Z^q=X}\xi(Z)$. This subsection was influenced by the
spirit in \cite{Dwork3} (in particular Section 3.5). This
subsection aims to prove Theorem \ref{T:Up} with the following
lemma.

\begin{lemma}\label{L:Dwork1}
Let $X\in\Omega_a$, $\hP\in\cO_a$, and $s
\in p^{\Z_{>0}}$.
If $r\geq  p^{-\frac{p}{p-1}}$
then \\
(1)
$|X-\hP|_p>r^{\frac{1}{s}}$ implies $|X^s-\hP^s|_p=|X-\hP|_p^s >r$,
and
(2) $|X^s-\hP^s|_p>r$ implies $|X^s-\hP^s|_p=|X-\hP|_p^s$.
\end{lemma}
\begin{proof}
(1)
If $\hP=0$ then the lemma is trivial.
We assume $\hP\neq 0$ for the rest of the proof.
Write $s=p^k$ for $k\geq 1$. We shall use induction on $k$.
We first prove the case $k=1$ for both statements.
That is,
$|X^p-\hP^p|_p>p^{-\frac{p}{p-1}}$ implies that $|X-\hP|_p>p^{-\frac{1}{p-1}}$
which in turn implies that $|X^p-\hP^p|_p=|X-\hP|_p^p$.

Write $Y:=X-\hP$. Then $X^p-\hP^p=(Y+\hP)^p-\hP^p= Y^p+pG$ where
$G= \sum_{m=1}^{p-1}(\binom{p}{m}/p)Y^{p-m}\hP^m
\in Y \Z[\hP,Y]$.
If $|Y|_p\leq p^{-\frac{1}{p-1}}$, that is,
$\ord_p Y\geq \frac{1}{p-1}$,
then $|pG|_p\leq |pY|_p\leq
p^{-\frac{p}{p-1}}$.
Since $|\hP|_p\leq 1$, one has
$\ord_p G \geq \ord_p Y$.
Thus $|X^p-\hP^p|_p\leq \max(|Y^p|_p,|pG|_p)
\leq p^{-\frac{p}{p-1}}$. Contradiction, so we have
$|Y|_p>p^{-\frac{1}{p-1}}$. This implies that
$\ord_p Y^{-i} > - \frac{i}{p-1}$ for any $i\in\Z$.
By the triangle inequality
$\ord_p G/Y^p \geq \min_{1\leq i\leq p-1} \ord_pY^{-i} > -1$.
Hence $\ord_p (pG) > \ord_p Y^p$.
Again by the triangle inequality
we have $|X^p-\hP^p|_p=|Y^p|_p=|X-\hP|_p^p$.

(2)  Suppose it holds for $s=p^{k-1}$. By assumption,
$|X^{p^k}-\hP^{p^k}|_p =
|(X^p)^{p^{k-1}}-(\hP^p)^{p^{k-1}}|_p>p^{-\frac{p}{p-1}}$. By
inductive argument one has $|X^{p^k}-\hP^{p^k}|_p=|X^p-\hP^p|_p^{p^{k-1}}$
and so $|X^p-\hP^p|_p>p^{-\frac{p}{p^{k-1}(p-1)}}\geq
p^{-\frac{p}{p-1}}$. The latter implies that $|X^p-\hP^p|_p=|X-\hP|_p^p$,
again by induction. Therefore, one has
$|X^{p^k}-\hP^{p^k}|_p=|X-\hP|_p^{p^k}$, as we desire. Now suppose
$|X-\hP|_p>r^{\frac{1}{p^k}}$. Then $|X-\hP|_p>p^{-\frac{1}{p-1}}$ and so
$|X^p-\hP^p|_p=|X-\hP|_p^p>r^{\frac{1}{p^{k-1}}}$. Then we use induction
argument to get
\begin{equation*}
|X^{p^k}-\hP^{p^k}|_p
=
|(X^p)^{p^{k-1}}-(\hP^p)^{p^{k-1}}|_p
=|X^p-\hP^p|_p^{p^{k-1}}
=|X-\hP|_p^{p^{k}}.
\end{equation*}
This finishes our proof.
\end{proof}

\begin{theorem}\label{T:Up}
Let $r\in |\Omega_a|_p$ and
$p^{-\frac{1}{p^{a-1}(p-1)}}<r<1$.
Let $s\in p^{\Z_{\geq 0}}$.
Then $U_p \cH(\bA_{r,s})\subseteq \cH(\bA_{r^p,sp})$.
Then $U_q=U_p^a$ and
$U_q\cH(\bA_{r})\subseteq \cH(\bA_{r^q})$.
\end{theorem}
\begin{proof}
1) We shall demonstrate a proof for the case
$s=1$ since the general case is very similar.
Let $\xi\in \cH(\bA_r)$.

Firstly, we show that $U_p\xi$ defines a function on
the affinoid $\bA_{r^p,p}$.
It suffices to show that
$Z^p=X\in \bA_{r^p,p}$ implies that $Z\in \bA_r$.
Indeed, for every $2\leq j\leq \ell$ one has $|Z^p-\hP_j^p|_p \geq r^p
> p^{-\frac{p}{p-1}}$ by hypothesis.
By Lemma \ref{L:Dwork1} one has
$|Z-\hP_j|_p^p=|Z^p-\hP_j^p|_p\geq r^p$. That is, $|Z-\hP_j|_p\geq r$.
On the other hand, by $|Z^p|_p \leq 1/r^{p}$,
one has $|Z|_p\leq 1/r$. This proves our claim.

Secondly we show that $U_p\xi\in\cH(\bA_{r^p,p})$.
Our proof below follows \cite[Lemma on page 40]{Dwork3}.
Before we start, an easy fact is prepared:
\begin{eqnarray}\label{E:norm}
\sup_{X\in \bA_{r^p,p}} |(U_p\xi)(X)|_p &\leq  &
p\cdot \sup_{X\in \bA_{r}}|\xi(X)|_p.
\end{eqnarray}
Let $\Tr$ denote the trace map from $\Omega_a(Z)$ to $\Omega_a(X)$
where $Z$ is a function with $Z^p=X$.
If $\xi\in \Omega_a(X)$, then by definition $U_p\xi = \frac{1}{p}\cdot
\Tr\;\xi(X)$. This shows that $U_p$ maps $\Omega_a(X)$ to
itself and by (\ref{E:norm}), if $\xi$ has no pole in $\bA_r$ then
$ U_p\xi$ has no pole in $\bA_{r^p,p}$. Thus
$U_p$ restricts to a mapping
$
\Omega_a(X)\cap \cH(\bA_r) \longrightarrow \Omega_a(X)\cap \cH(\bA_{r^p,p}),
$
which is continuous relative to the supremum norms.
Since $\xi\in \cH(\bA_r)$, one gets that $\xi$
may be uniformly approximated on $\bA_r$ by elements
of $\Omega_a(X)\cap\cH(\bA_r)$ and so by (\ref{E:norm}) again
$U_p\xi$ can be uniformly approximated on $\bA_{r^p,p}$
by elements
of $\Omega_a(X)\cap \cH(\bA_{r^p,p})$.
This completes the proof of the assertion about $U_p$.

2)
Let $Z^q=X\in \bA_{r^q}$ for $r>p^{-\frac{1}{p^{a-1}(p-1)}}$.
For $2\leq
j\leq \ell$, one has $|Z^q-\hP_j^q|_p\geq r^q$; so
$|Z-\hP_j|_p=|Z^q-\hP_j^q|_p^{1/q}\geq r$. One also
observes that $|Z^q|_p\leq 1/r^q$ implies that $|Z|_p\leq 1/r$. This
proves that $Z\in \bA_r$. This proves that
$U_q\xi$ is a function on $\bA_{r^q,q}$.
As $\hP_j^q=\hP_j$ for all $j$ one has $\bA_{r^{q},q}=\bA_{r^{q}}$,
it follows that $U_q\xi$ is defined over $\bA_{r^q}$.
\end{proof}

\subsection{Push-forward maps and Dwork's splitting functions}
\label{S:2.3}

In the previous subsection
we have defined the $U_p$ and $U_q$ operators on
suitable $p$-adic Banach spaces.
It remains to define the ``Dwork's splitting function" to
finish the process of defining the Frobenius map.
Let $\tau$ be a lift of the Frobenius endomorphism $c\mapsto c^p$
of $\F_{p^a}$ to $\Omega_a$ which fixes $\Omega_1$.
Thus $\tau$ generates $\Gal(\Omega_a/\Omega_1)$.
Let $\bA_{r,s}^\tau$ denote the image of $\bA_{r,s}$ under $\tau$.

Since one may have a pole $\hP_j$ other than $0$ and $\infty$,
one encounters the following problem:
for any $\xi(X)\in \cH(\bA_{r})$
its image $\xi^\tau(X)$ does not lie in
$\cH(\bA_{r}^\tau)$ anymore.
So the naive generalization of Dwork's splitting function does not work.
This prompts us to define some push-forward maps.

Define a map of $p$-adic Banach spaces
\begin{eqnarray*}
\cH(\bA_{r,s}) &\stackrel{\tau_*}{\longrightarrow} & \cH(\bA_{r,s}^\tau)\\
\xi          &\mapsto &\tau\circ \xi\circ \tau^{-1}.
\end{eqnarray*}
For any $k\in \Z_{\geq 0}$, one has
$\tau_*(\cH(\bA_{r,p^k}))=\cH(A^\tau_{r,p^k})=\cH(\bA_{r,p^{k+1}})$.
As a simple example, for $B\in \Omega_a$ and a Teichm\"uller lift
$\hP$ in $\Omega_a$ with $\xi(X)=\frac{B}{X-\hP}\in \cH(\bA_r)$ we
have $(\tau_*\xi)(X)=
\frac{\tau(B)}{X-\tau(\hP)}=\frac{\tau(B)}{X-\hP^p}$. On the other
hand, one may check routinely that $\tau_*^k$ commutes with $U_p$
for any $k\in\Z$.

For any $\hat{f}(x)$ (fixed
in Section \ref{S:2.0}),
and for every $1\leq j\leq \ell$,
let
\begin{eqnarray}\label{E:F_j}
F_j(X_j)&:=& \prod_{i=1}^{d_j}E(\gamma a_{j,i}X_j^i)
\end{eqnarray}
where we recall that $E(X)$ is the Artin-Hasse
exponential function
and $\gamma$ is the root of $\log E(X)$ with
$\ord_p \gamma = \frac{1}{p-1}$.
We now induce our new splitting functions:
\begin{equation}\label{E:F}
F(X):= \prod_{j=1}^{\ell}F_j(X_j); \quad \quad
F_{[a]}(X):= \prod_{k=0}^{a-1}(\tau_*^k F)(X^{p^k}).
\end{equation}

\begin{lemma}\label{L:FX}
Let $k\geq 1$ be any integer. Let $r\in |\Omega_a|_p$ and
$p^{-\frac{p}{p-1}}<r<1$. Then for any $\xi(X)\in
\cH(\bA^{\tau^k}_{r})$ one has $\xi(X^{p^k})\in
\cH(\bA_{r^{1/p^k}})$.
\end{lemma}
\begin{proof}
It suffices to show that
$|X-\hP_j|_p\geq r^{1/p^k}$ implies that $|X^{p^k}-\hP_j^{p^k}|_p\geq r$.
This follows from Lemma \ref{L:Dwork1}
immediately.
\end{proof}

\begin{theorem}\label{T:FX}
Let $d_0:=\max_{1\leq j\leq \ell}d_j$.
Let $r\in |\Omega_a|_p$ and
$p^{-\frac{1}{d_0p^{a-1}(p-1)}} <  r < 1$. Then $F(X)\in
\cH(\bA_{r^{p^{a-1}}})\subseteq\cH(\bA_r)$ and $F_{[a]}(X)\in \cH(\bA_r)$.
\end{theorem}
\begin{proof}
Write $r_j:=p^{-\frac{1}{d_j(p-1)}}$.
Write $F_j(X)=\sum_{n=0}^{\infty}F_{j,n}X^n$ over $\cO_a$.
Note that $F_j(X)$'s convergence
radius is
$\liminf_{n}|F_{j,n}|_p^{-1/n}=
p^{\liminf_{n} \ord_pF_{j,n}/n} \geq 1/r_j.
$
(See Lemma \ref{L:F}.)
By hypothesis, one has
$r^{p^{a-1}}> r_j$ for every $j$,
so $F_1(X), F_j((X-\hP_j)^{-1})\in\cH(\bA_{r^{p^{a-1}}})$
for $2\leq j\leq \ell$.
Hence $F(X)\in\cH(\bA_{r^{p^{a-1}}})$.

Then $(\tau_*^k F)(X)\in\cH(A^{\tau^k}_{r^{p^{a-1}}})$ for every
$0\leq k\leq a-1$. Our hypothesis implies that
$p^{-\frac{p}{p-1}}<r^{p^{a-1}}<1$.
So one may apply  Lemma \ref{L:FX} and gets,
$(\tau_*^k F)(X^{p^k})\in\cH(\bA_{r^{p^{a-1-k}}})\subseteq \cH(\bA_r)$
for every $0\leq k\leq a-1$.
Therefore, their product $F_{[a]}(X)$ lies in $\cH(\bA_r)$ as well.
\end{proof}

\subsection{The trace formula of $\alpha_a$}\label{S:2.4}

For the rest of the paper we assume
\begin{equation}\label{E:r}
\fbox{$r\in |\Omega_a|_p$ and $p^{-\frac{1}{d_0p^{a-1}(p-1)}} <  r  < 1$
where $d_0=\max_{1\leq j\leq \ell} d_j$.}
\end{equation}
This bound of $r$ is to assure that $F_{[a]}(X)$ lies in $\cH(\bA_r)$.
Let $\alpha_a:=U_q\circ F_{[a]}(X)$, by which we mean
the composition map of
$U_q$ with the multiplication map by $F_{[a]}(X)$. Then $\alpha_a$ is a
$\Omega_a$-linear map from $\cH(\bA_r)$ to $\cH(\bA_{r^q})$ by
Theorem \ref{T:Up}. Composing with the natural restriction map
$\cH(\bA_{r^q})\rightarrow\cH(\bA_r)$, one observes that
$\alpha_a$ defines an endomorphism of $\cH(\bA_r)$.

\begin{lemma}[Dwork-Monsky-Reich]\label{L:Monsky-Reich}
Let $\bar{f}\in\A(\F_q)$.
Let $r$ be as in (\ref{E:r}), then  the $\Omega_a$-linear
endomorphism $\alpha_a$ of $\cH(\bA_r)$ is completely continuous
and one has
\begin{eqnarray}\label{E:L}
L(\bar{f};T)&=&\frac{\det(1-T\alpha_a| \cH(\bA_r ))} {\det(1-T
q\alpha_a|\cH(\bA_r ))}.
\end{eqnarray}
\end{lemma}
\begin{proof}
Let $\cH^\dagger(\bA_1):=\bigcup_{0<r<1}\cH(\bA_r)$.
One notes that
$\cH^\dagger(\bA_1)$ is the Monsky-Washnitzer dagger space of $\bA_1$.
Our assertions then follow from
the trace formula of \cite{Monsky-Washnitzer:I} and \cite{Reich},
as explained
in \cite[Section 6]{Robba:1} and see (6.3.11) in particular.
(Our hypothesis $g(0)=0$ was used there).
Basically their trace formula says that $\alpha_a$
is completely continuous on $\cH^\dagger(\bA_1)$ and
$\det(1-T
\alpha_a|\cH^\dagger(\bA_1)) = \det(1-T \alpha_a|\cH(\bA_r))$
for any $r$ within our range in (\ref{E:r}).
Since it is routine to check this, we omit details.
\end{proof}

\begin{remark}
One can also formulate the above trace formula using
Berthelot's rigid cohomology theory.
See \cite{Berthelot} for detailed annotation of
Robba's formulation in \cite[Section 6]{Robba:1}.
\end{remark}

\subsection{Descent from $\alpha_a$ to $\alpha_1$}\label{S:2.5}

The results in this subsection are only used in the proof of
Theorem \ref{T:main1} in Section \ref{S:4}. Below we shall use a
subindex in $\det_{\Omega_1}(\cdot)$ or $\det_{\Omega_a}(\cdot)$
to emphasize our consideration of a map over $\Omega_1$-space or
$\Omega_a$-space, respectively. We shall omit the base space
$\cH(\bA_r)$ in $\det(\cdot)$ if the context clearly assures that
no confusion is possible. The upshot of our argument is to
``descent" the $\alpha_a$ map of the $\Omega_a$-space $\cH(\bA_r)$
to the $\alpha_1$ map of the $\Omega_1$-space $\cH(\bA_r)$. This
idea appeared initially in \cite[Section 7]{Dwork2}. In this paper
we use $\NP_p(\cdot)$ and $\NP_q(\cdot)$ to denote $p$-adic and
$q$-adic Newton polygons, respectively. (These should not be
confused with $\NP(\bar{f};\F_q)$.) We use $1/a\cdot\NP_p(\cdot)$
to denote the image of $\NP_p(\cdot)$ shrunk by a factor of $1/a$.

\begin{lemma}\label{L:alpha}
Let $ \alpha_1:=\tau_*^{-1}\circ U_p\circ F(X)$. Then $\alpha_1$
is a completely continuous $\tau^{-1}$-linear map from
$\cH(\bA_r)$ to $\cH(\bA_{r^p})$ (over $\Omega_a$) and $\alpha_a =
\alpha_1^a$ as $\Omega_1$-linear maps. Then
\begin{eqnarray}\label{E:det}
\det_{\Omega_a}(1-T^a\alpha_a)^a &=&
\prod_{k=0}^{a-1}\det_{\Omega_1} (1-T\zeta_a^k\alpha_1),
\end{eqnarray}
where $\zeta_a$ is a primitive $a$-th root of unity. Then $
\NP_q(\det_{\Omega_a}(1-T\alpha_a))=1/a\cdot
\NP_p(\det_{\Omega_1}(1-T\alpha_1)). $
\end{lemma}

\begin{proof}
As we already remarked at the beginning of Section \ref{S:2.3},
$\tau_*^{-1}$ and $U_p$
commute with each other.
For any $k\in\Z$, the
$\Omega_1$-linear multiplication map of $(\tau_*^k F)(X)$
on $\cH(A^{\tau^k}_r)$
can be written as
$\tau_*^k\circ F(X)\circ \tau_*^{-k}$.
On the other hand,
for any function $H_k(X)\in \cH(\bA_r^{\tau^k})$ one has a general
identity stating that
\begin{eqnarray}\label{E:Bombieri}
U_q\circ\prod_{k=0}^{a-1}H_k(X^{p^k})
&=& \prod_{k=0}^{a-1} U_p\circ H_{a-1-k}(X),
\end{eqnarray}
where second product is noncommutative and its factors are
ordered from left to right as $k$ increases. We retain this notation of
noncommutative products for the rest of the paper.

Now apply (\ref{E:Bombieri})
to $F_{[a]}(X)$ with $H_k(X):=(\tau_*^k F)(X)$.
One has
\begin{eqnarray*}
U_q\circ F_{[a]}(X)
&=&
    \prod_{k=0}^{a-1} (U_p \circ \tau_*^{a-1-k}\circ F(X)\circ \tau_*^{-(a-1-k)})\\
&=& \prod_{k=0}^{a-1} (\tau_*^{a-1-k}\circ U_p\circ F(X)\circ \tau_*^{-(a-1-k)}).
\end{eqnarray*}
By telescoping, one  gets $U_q\circ F_{[a]}(X) = (\tau_*^{-1}\circ
U_p\circ F(X))^a$. That is, $\alpha_a = \alpha_1^a$.

The proof for $\alpha_1$ being completely continuous is verbatim
for $\alpha_a$ which is already proved in Lemma
\ref{L:Monsky-Reich}. Now it is elementary to see that
$$
\det_{\Omega_1}(1-T^a\alpha_1^a) =
\prod_{k=0}^{a-1}\det_{\Omega_1}(1-T\zeta_a^k\alpha_1).
$$
One may also show as an exercise that (see \cite[(41)]{Bombieri} for
details)
$$
\det_{\Omega_a}(1-T\alpha_a)^a=\det_{\Omega_1}(1-T\alpha_a).
$$
Combining these two equalities with $\alpha_a=\alpha_1^a$, one
obtains (\ref{E:det}). The last assertion about Newton polygons
follows from the elementary theory of Newton polygons (see
\cite[Lemma 1.6]{Dwork1} and \cite[Lemma 7.1]{Dwork2}).
\end{proof}

\begin{proposition}\label{P:reduce1}
The slope $<1$ part (of horizontal length $d-\ell+1$) of
$\NP(\bar{f};\F_q)$ is equal to
$\NP_q(\det_{\Omega_a}(1-T\alpha_a)\bmod T^{d-\ell+1})$ which is
equal to
$$1/a\cdot
\NP_p(\det_{\Omega_1}(1-T\alpha_1)\bmod T^{a(d-\ell+1)+1}).
$$
The same holds if one replaces $\Omega_1$ and $\Omega_a$
by $\Omega_1'$ and $\Omega_a'$, respectively.
\end{proposition}

\begin{proof}
By (\ref{E:L}), one has
\begin{eqnarray}\label{E:proposition}
L(f\bmod\cP;T)\cdot\det_{\Omega_a}(1-Tq\alpha_a)&=&\det_{\Omega_a}(1-T\alpha_a).
\end{eqnarray}
Note that all slopes are greater than or equal to $1$ in
$\NP_q(\det_{\Omega_a}(1-Tq\alpha_a))$. By the Weil conjectures
for (projective) curves, $L(\overline{f};F_q)$ is a degree $d$ polynomial
(see (\ref{E:20})) with all slopes in $[0,1]$. The slope-1 part of
$\NP(\bar{f};\F_q)$ is precisely of horizontal length $\ell-1$
(see Remark \ref{R:main2}). Let $\lambda$ be the biggest slope of
$\NP(\bar{f};\F_q)$ that is strictly less than $1$. Then the slope
$\leq \lambda$ part of $\NP(\bar{f};\F_q)$ is equal to
$\NP_q(\det_{\Omega_a}(1-T\alpha_a)\bmod T^{d-\ell+2})$ by
(\ref{E:proposition}) and the $p$-adic Weierstrass preparation
theorem (see section \cite[IV.4]{Koblitz}).

By Lemma \ref{L:alpha}, $1/a\cdot
\NP_p(\det_{\Omega_1}(1-T\alpha_1))
=\NP_q(\det_{\Omega_a}(1-T\alpha_a))$. By the previous paragraph,
the latter polygon has a vertex point at $T^{d-\ell+1}$, which
separates the slope $\leq \lambda$ and slope-$1$ segments. Hence
the former polygon has a corresponding vertex point at
$T^{a(d-\ell+1)}$. The upshot is that
$$1/a\cdot \NP_p(\det_{\Omega_1}(1-T\alpha_1)\bmod T^{a(d-\ell+1)+1})
=\NP_q(\det_{\Omega_a}(1-T\alpha_a)\bmod T^{d-\ell+2}).$$
Compiling these two paragraphs, our assertion follows. The last
assertion is obvious.
\end{proof}

\section{$p$-adic estimates of $L$-functions of exponential sums}
\label{S:3}

This section aims to prove Theorem \ref{T:reduce1} whose proof
is however very technical, so the reader
is recommended to refer to it
only when needed.
We retain all notations from previous sections, in particular
we recall the two bases
$\vec{b}_{\rm unw}$ and $\vec{b}_{\rm w}$
of $\cH(\bA_r)$ from Lemma \ref{L:basis}.
For any $c\in\R$ we denote by $\lceil c\rceil$ the least integer
greater than or equal to $c$.

We start with a lemma
inspired by a ``Dwork's Lemma" in \cite{Dwork3}:

\begin{lemma}\label{L:Dwork2}
Let $m\geq 1$ and $J\geq 3$. Let $r$ be as in (\ref{E:r})  of
Section \ref{S:2.4}. Then for any $X\in \bA_r$ one has
\begin{eqnarray}\label{E:Dwork2}
U_p(X-\hP_J)^{-m} &=& \sum_{n=\lceil{m/p}\rceil}^{m}
C^{n,m}\hP_J^{np-m}(X-\hP_J^p)^{-n},
\end{eqnarray}
where $C^{n,m}\in \Z_p$ with $\ord_p C^{n,m} \geq
\frac{np-m}{p-1}-1$. Then $\ord_p C^{n,m} = 0$ if and only if
$n=\lceil{\frac{m}{p}}\rceil$. If this is the case, then one has
 $C^{n,m}\equiv (-1)^{\epsilon-1}
\bmod p$ where $\epsilon=m-(n-1)p$.
\end{lemma}
\begin{proof}
By Theorem \ref{T:Up} one has $U_p(X-\hP_J)^{-m}\in\cH(\bA_{r^p,p})$.
The first statement of this lemma follows from an analogous
verification as presented in Section 5.3 of \cite{Dwork3}, so we
shall omit its proof here. We shall prove our last assertion
below.

By \cite[Lemma of Section 5.3, page 74]{Dwork3},
one has $C^{n,m} = \frac{m}{np}
\sum_{\vec{i}}\prod_{k=1}^{n}\binom{p}{i_k}$ where
$\vec{i}:=(i_1,\ldots, i_n)$ ranges in $\Z^n$ with $1\leq
i_1,\ldots,i_n\leq p$ and $\sum_{k=1}^{n}i_k = m$ (we denote this
set of $\vec{i}$ by $\cI$). Write $m=(n-1)p+r$ for some $1\leq
r\leq p$. If $r=p$ then one can easily see that $C^{n,m}=1$ and
our assertion clearly holds. Below we let $1\leq r\leq p-1$. We
assume additionally that $p>2$ since if $p=2$ then one has
$m=2n-1$ and it is easy to check that $\ord_pC^{n,m}=0$ directly.
Write $\varsigma(\vec{i}):=\prod_{k=1}^{n}\binom{p}{i_k}$. For any
$1\leq t\leq n$ let $\cI_t$ be the subset of $\cI$ consisting of
all $\vec{i}$ with $\ord_p(\varsigma(\vec{i}))=t$. It is clear
that $\cI=\bigcup_{t=1}^{n}\cI_t$ is a partition of $\cI$. Since
$\ord_p\binom{p}{i_k}=0$ (resp., $=1$) if and only if $i_k=p$
(resp., $\neq p$), one gets for any $1\leq t\leq n$ that
$\vec{i}\in \cI_t$ if and only if the $\vec{i}$ contain precisely
$t$ non-$p$ components. For each $\vec{i}\in\cI_t$ there are
actually $\binom{n}{t}$ of them by a permutation of the non-$p$
components among the $n$ components and they have the same
$\varsigma(\vec{i})$. Let $\cJ_t$ be the set of all $t$-tuples
$\vec{j}:=(j_1,\ldots,j_t)$ with $1\leq j_1,\ldots,j_t\leq p-1$
and $\sum_{k=1}^{t}j_k = (t-1)p+r$. Then one gets
\begin{eqnarray}\label{E:Cnm}
C^{n,m} &=& \sum_{t=1}^{n}\sum_{\vec{i}\in\cI_t}\frac{m}{np}\varsigma(\vec{i})
 = \sum_{t=1}^{n}\sum_{\vec{j}\in\cJ_t}
     \frac{m}{np}\binom{n}{t}\varsigma(\vec{j})\\\nonumber
 &=& \frac{m\binom{p}{r}}{p} +
     \sum_{t=2}^{n}\sum_{\vec{j}\in\cJ_t}
     \frac{m}{np}\binom{n}{t}
     \varsigma(\vec{j}).
\end{eqnarray}
One easily observes that the first summand
is a $p$-adic unit.
Now we claim that for any $t\geq 2$ and $\vec{j}\in \cJ_t$ one
has $\ord_p(\frac{m}{np}\binom{n}{t}
     \varsigma(\vec{j}))\geq 1$.
Indeed, one has for some $u\in Z_p$ that
$$
\frac{m}{np}\binom{n}{t}\varsigma(\vec{j})
= \frac{m}{np}\left(\frac{n}{t}\binom{n-1}{t-1}\right) ( p^t u)
= um\binom{n-1}{t-1}\frac{p^{t-1}}{t}.
$$
It is easy to observe that for any $t\geq 2$ and $p>2$ one has
$\ord_p t\leq t-2$. This proves our claim above.
By (\ref{E:Cnm}), we get that $C^{n,m}$ is a $p$-adic unit.
\end{proof}

The computation of $\alpha_1=\tau_*^{-1}\circ U_p\circ F(X)$ uses
the observation that $\tau_*^{-1}$ and $U_p$ respect the
Mittag-Leffler decomposition while the multiplication map $F(X)$
does not. For $1\leq j\leq \ell$ and for any
$\xi(X)\in\cH(\bA_r)$, let $\xi(X)_{\hP_j}$ denote the $j$-th
component in the $p$-adic Mittag-Leffler decomposition as in Lemma
\ref{L:KML}. We recall our notation $X_1=X$ and
$X_j=(X-\hP_j)^{-1}$ for $2\leq j\leq \ell$.

Now we recall certain properties of $F_j(X_j) =
\sum_{n=0}^{\infty}F_{j,n}X_j^n\in \cO_a[[X_j]]$ (see proof in,
for instance, \cite[Section 1]{Zhu:1}).

\begin{lemma}\label{L:F}
For any $1\leq j\leq \ell$ and $n\geq 0$ one has
$\ord_p F_{j,n}
 \geq  \frac{\left\lceil\frac{n}{d_j}\right\rceil}{p-1}
$
where the equality holds if $d_j|n$ and
$\frac{n}{d_j}\leq p-1$.
In particular, $\ord_pF_{j,n}>0$ for any $n>0$.
\end{lemma}

\begin{lemma}[Key computational lemma]\label{L:computation}
(1) If $\xi(X)\in\cH(\bA_r)$ is given by its Laurent expansion at
$\hP_J$, that is $\xi(X)=\sum_{n=-\infty}^{\infty}B_nX_{J}^n$ for
some $B_n\in\Omega_a$, then
$(\xi(X))_{\hP_{J}}=\sum_{n=0}^{\infty}B_nX_{J}^n$, and $B_0=0$ if
$\hP_{J}\neq \infty$.

(2)
Recall $F(X)$ and $C^{n,m}$ from (\ref{E:F}) and (\ref{E:Dwork2})
respectively.
For all $i\geq 0$ write
$(F(X)X_J^i)_{\hP_{J_1}} = \sum_{n=0}^{\infty}H_{J_1,J}^{n,i}
                          X_{J_1}^n$
for some $H_{J_1,J}^{n,i}\in \Omega_a$.
Then
$$(\alpha_1 X_J^i)_{\hP_{J_1}}
=\sum_{n=0}^{\infty}B_{J_1,J}^{n,i}
X_{J_1}^n \in\Omega_a[[X_{J_1}]]
$$
where
$$
B_{J_1,J}^{n,i}:=
\left\{
\begin{array}{ll}
 \tau^{-1}H_{J_1,J}^{np,i}
     &\mbox{for $J_1=1,2$}    \\
 \sum_{m=n}^{np}C^{n,m}\hP_{J_1}^{n-mp^{a-1}}
      \tau^{-1}H_{J_1,J}^{m,i}
  &\mbox{for $J_1\geq 3$.}
\end{array}
\right.
$$
\end{lemma}
\begin{proof}
Part (1) is a simple corollary of the remarks preceding the lemma and
Lemmas \ref{L:KML} and \ref{L:Dwork2}. The rest are routine consequences.
\end{proof}

For any integers $s\geq 0$ and $t\geq 1$ we use $\cC(s,t)$ to
denote the condition that $t|s$ and $0\leq \frac{s}{t}\leq p-1$
are satisfied (e.g., the condition in Lemma \ref{L:F} is
$\cC(n,d_j)$). We claim the following:
\begin{equation}
    \label{E:rough}
|H_{J_1,J}^{n,i}|_p \leq
\left\{
    \begin{array}{ll}
    p^{-\frac{n-i}{d_{J_1}(p-1)}}      &\mbox{ if $J_1=J$}\\
    p^{-\frac{n+i}{d_{J_1}(p-1)}}      &\mbox{ if $J_1=1\neq J$}\\
    p^{-\frac{n}{d_{J_1}(p-1)}}        &\mbox{ if $J_1=2\neq J$.}
    \end{array}
\right.
\end{equation}
Furthermore, the equalities hold
if and only if additional conditions $\cC(n-i,d_{J_1})$,
$\cC(n+i,d_{J_1})$, $\cC(n,d_{J_1})$ hold, respectively.

A proof for the case $J=J_1$ is sketched below and proofs for
other cases are omitted as they are formal and similar. Let
$\vec{n}:=(n_1,\ldots,n_\ell)\in\Z_{\geq 0}^\ell$. Then  one
notices that for $J=1$
\begin{eqnarray*}
H_{1,1}^{n,i}&=& \sum \left( F_{1,n_1}\prod_{j\neq 1}
  \left(\sum_{m_j=0}^{n_j} F_{j,m_j}\binom{n_j-1}{m_j-1}\hP_j^{n_j-m_j}
  \right)\right),
\end{eqnarray*}
and for $J\geq 2$,
\begin{eqnarray*}
H_{J,J}^{n,i}&=& \sum \left(
F_{J,n_J}\left(\sum_{m_1=n_1}^{\infty}F_{1,m_1}\binom{m_1}{n_1}
                 \hP_{J}^{m_1-n_1}\right)\right.\\
&&\left.\cdot\prod_{j\neq 1,J} \left(\sum_{m_j=0}^{\infty}
F_{j,m_j}(-1)^{m_j}
  \binom{n_j+m_j-1}{m_j-1}(\hP_j-\hP_{J})^{-(n_j+m_j)}
  \right)
  \right)
\end{eqnarray*}
where the sums both range over all $\vec{n}\in\Z_{\geq 0}^\ell$
such that $n-i=n_{J}-\sum_{j\neq J}n_j$.

{}From the above, one observes that $\ord_pH_{J_1,J}^{n,i}$ is
greater than or equal to the minimal valuation among the
$\vec{n}$-summand in its formula as $\vec{n}$ varies in its
domain. Each $\vec{n}$-summand is the product of $\ell$ elements
in $\cO_a$, so its valuation is equal to the sum of the valuations
of these $\ell$ elements in $\cO_a$. It is easy to observe that
$\ord_pH_{J_1,J}^{n,i}\geq \min_{n_1}(\ord_pF_{1,n_1}) \geq
\frac{n-i}{d_{J_1}(p-1)}$ as the minimum is taken over all
$n_1=n-i+\sum_{j=2}^{\ell}n_j$. Moreover, if $\cC(n-i,d_{J_1})$
holds then by Lemma \ref{L:F} the minimal is uniquely achieved at
$\vec{n}=(n-i,0,\ldots,0)$ and the equality holds. Conversely,
suppose the equality in (\ref{E:rough}) holds. It can be easily
seen that $H_{J_1,J}^{n,i}$ lies in $\cO_a$, in which $p$ has
ramification index $p-1$ over $\Z$, so $\ord_p
H_{J_1,J}^{n,i}=\frac{n-i}{d_{J_1}(p-1)}$ lies in
$\frac{\Z}{p-1}$. Thus $\cC(n-i,d_{J_1})$ holds. This proves our
claim in (\ref{E:rough}).

\begin{theorem}[Unweighted estimates]\label{T:sharp}
Let $B_{J_1,J}^{n,i}\in\cO_a$ be as in Lemma \ref{L:computation}(2).

(1) For $J_1=1,2$ and for $n,i\geq 0$ one has
$$
|B_{J_1,J}^{n,i}|_p \leq
\left\{
\begin{array}{ll}
p^{-\frac{np-i}{d_{J_1}(p-1)}} & \mbox{if $J_1=J$} \\
p^{-\frac{np+i}{d_{J_1}(p-1)}} & \mbox{if $J_1=1\neq J$}\\
p^{-\frac{np}{d_{J_1}(p-1)}} & \mbox{if $J_1=2\neq J$}.
\end{array}
\right.
$$
The equalities hold if and only if
the additional conditions
$\cC(np-i,d_{J_1}), \cC(np+i,d_{J_1})$, and $\cC(np,d_{J_1})$
hold, respectively.

(2) For $J_1\geq 3$ and for $n\geq 1, i\geq 0$ one has
$$
|B_{J_1,J}^{n,i}|_p
\leq
\left\{
\begin{array}{ll}
p^{-\frac{(n-1)p-(i-1)}{d_{J_1}(p-1)}}
                            & \mbox{if $J_1=J$} \\
p^{-\frac{(n-1)p+1}{d_{J_1}(p-1)}}
                            & \mbox{if $J_1\neq J$}.
\end{array}
\right.
$$
For $d_{J_1}\geq 2$ the equalities hold
if additional conditions
$\cC((n-1)p-(i-1),d_{J_1})$
and $\cC((n-1)p+1,d_{J_1})$ hold, respectively.
\end{theorem}

\begin{proof}
If $J_1=1,2$ one has
$|B_{J_1,J}^{n,i}|_p = |H_{J_1,J}^{np,i}|_p$
by Lemma \ref{L:computation}. Combining this with (\ref{E:rough}),
part (1) follows immediately. We are left to prove part (2).
Assume $J_1\geq 3$ from now on. We shall outline a proof for the
case $J=J_1$: Let $n,i$ be fixed in their appropriate ranges. By
Lemma \ref{L:computation}(2), one has that
$|B_{J_1,J}^{n,i}|_p\leq \max_{n\leq m\leq np}
(|H_{J_1,J}^{m,i}C^{n,m}|_p) $ and the equality holds if the
maximum is unique. Pick $m_0:=(n-1)p+1$, then one has two cases:
(a) For any $m_0< m\leq np$
one has
$|H_{J_1,J}^{m,i}|_p < p^{-\frac{m_0-i}{d_{J_1}(p-1)}}$
by (\ref{E:rough});
(b) Let $n\leq m\leq m_0$. The function
$c(m):=\frac{m-i}{d_{J_1}(p-1)}+(\frac{np-m}{p-1}-1)
=\frac{np}{p-1}-\frac{i}{d_{J_1}(p-1)}-1
- m\frac{d_{J_1}-1}{d_{J_1}(p-1)}
$
has its minimum $c(m_0)$.
If $d_{J_1}\geq 2$ then this minimum is unique.
By (\ref{E:rough}) and Lemma \ref{L:Dwork2} one has
$$
\max_{n\leq m\leq m_0}|H_{J_1,J}^{m,i}C^{n,m}|_p
\leq p^{-c(m_0)}=p^{-\frac{m_0-i}{d_{J_1}(p-1)}}
=p^{-\frac{(n-1)p-(i-1)}{d_{J_1}(p-1)}}.
$$
Combining (a) and (b) one gets the desired upper bound
for $|B_{J_1,J}^{n,i}|_p$.
Suppose $d_{J_1}\geq 2$ and $\cC(m_0-i,d_{J_1})$ holds. Then
the maximum is achieved uniquely at $m_0$ by Lemma
\ref{L:Dwork2}.
Combining the above, we have proved part (2) for the case $J_1=J$.
Since other cases are similar and we omit them here, and
finally we conclude the proof to our theorem.
\end{proof}

\begin{theorem}[Weighted estimates]\label{T:reduce1}
Write $(\alpha_1 Z_J^i)_{\hP_{J_1}} = \sum_{n=0}^{\infty}
C_{J_1,J}^{n,i}Z_{J_1}^n $ in $\Omega'_a[[Z_{J_1}]]$. Then
$C_{J_1,J}^{n,i}
=B_{J_1,J}^{n,i}\gamma^{\frac{i}{d_J}-\frac{n}{d_{J_1}}}$ and $
\ord_pC_{J_1,J}^{n,i} = \ord_p B_{J_1,J}^{n,i} +
\frac{1}{p-1}(\frac{i}{d_J}-\frac{n}{d_{J_1}}) $ where
$B_{J_1,J}^{n,i}$ is as in Lemma \ref{L:computation} (2).

(1) For $J_1=1,2$ and $n,i\geq 0$ one has
\begin{equation}\label{E:bound12}
\ord_p C_{J_1,J}^{n,i} \geq
\left\{
\begin{array}{ll}
\frac{n}{d_{J_1}}
    &\mbox{ if $J_1=J$}\\
\frac{n}{d_{J_1}}+\frac{i}{p-1}(\frac{1}{d_{J_1}}+\frac{1}{d_J})
    &\mbox{ if $J_1=1\neq J$}\\
\frac{n}{d_{J_1}}+\frac{i}{(p-1)d_J}
    &\mbox{ if $J_1=2\neq J$}.
\end{array}
\right.
\end{equation}
The equalities hold if and only if
$\cC(np-i,d_{J_1}), \cC(np+i,d_{J_1})$ and $\cC(np,d_{J_1})$ hold,
respectively.

(2) For $J_1\geq 3$ and $n\geq 1, i\geq 0$ one has
\begin{equation}\label{E:bound3}
\ord_p C_{J_1,J}^{n,i} \geq
\left\{
\begin{array}{ll}
\frac{n-1}{d_{J_1}}
    &\mbox{ if $J_1=J$}\\
\frac{n-1}{d_{J_1}}+\frac{i}{(p-1)d_{J}}
    &\mbox{ if $J_1\neq J$}.
\end{array}
\right.
\end{equation}
For $d_{J_1}\geq 2$ the equalities hold if  conditions
$\cC((n-1)p-(i-1),d_{J_1}), \cC((n-1)p+1,d_{J_1})$ hold, respectively.
\end{theorem}

\begin{proof}
The first statement is clear by a simple calculation.
Parts (1) and (2) follow from Lemma \ref{L:F}
and parts (1) and (2) of Theorem \ref{T:sharp}, respectively.
\end{proof}

Notations: let $j,j_1\geq 3$, $j\neq j_1$, $n,i\geq 1$.
We put row minimal $p$-adic valuation in boxes.
\begin{table}[ht]
  \centering\label{Table:reduce1}
\begin{tabular}{|c||llll|}
\hline
$\ord_p(\cdot)\geq$
    & $Z_{J=1}^i$
    & $Z_{J=2}^i$
    & $Z_{J=j}^i$
    & $Z_{J=j_1}^i$
    \\\hline\hline
   $Z_{J_1=1}^n$
       & \fbox{$\frac{n}{d_1}$}
       & $\frac{n}{d_1}+ \frac{i/d_1+i/d_2}{p-1}$
       & $\frac{n}{d_1}+ \frac{i/d_1+i/d_{j}}{p-1}$
       & $\frac{n}{d_1}+ \frac{i/d_1+i/d_{j_1}}{p-1}$
       \\ 
   $Z_{J_1=2}^n$
       & $\frac{n}{d_2}+\frac{i/d_1}{p-1}$
       & \fbox{$\frac{n}{d_2}$}
       & $\frac{n}{d_2}+ \frac{i/d_{j}}{p-1}$
       & $\frac{n}{d_2}+ \frac{i/d_{j_1}}{p-1}$
       \\   
   $Z_{J_1=j}^n$
       & $\frac{n-1}{d_{j}}+\frac{i/d_1}{p-1}$
       & $\frac{n-1}{d_{j}}+\frac{i/d_2}{p-1}$
       & \fbox{$\frac{n-1}{d_{j}}$}
       & $\frac{n-1}{d_{j}}+\frac{i/d_{j_1}}{p-1}$
       \\
   $Z_{J_1=j_1}^n$
       & $\frac{n-1}{d_{j_1}}+\frac{i/d_1}{p-1}$
       & $\frac{n-1}{d_{j_1}}+\frac{i/d_2}{p-1}$
       & $\frac{n-1}{d_{j_1}}+\frac{i/d_{j}}{p-1}$
       & \fbox{$\frac{n-1}{d_{j_1}}$}
       \\   
  \hline
\end{tabular}
  \caption{Lower bounds for $\ord_pC_{J_1,J}^{n,i}$ in matrix of $\alpha_1$}
  \end{table}

\section{Newton polygon lies over Hodge polygon}\label{S:4}

Our proof of Theorem \ref{T:main1}
consists of three parts. The first two parts are in the
spirit of Dwork (see \cite[Section 7]{Dwork2} or \cite[Lemma
2]{Bombieri}) after a simple reduction. The third part uses Wan's
\cite[Theorem
2.4]{Wan:1}.

We queue up the numbers in (\ref{E:HP}) in nondecreasing order.
For any $i\geq 1$, let $m_i$ be the $i$-th in this queue.
For any $k\geq 1$, let $c_k:=\sum_{i=1}^{k}m_i$ and set $c_0=0$.
It is by elementary arithmetic of Newton polygons that
$\HP(\A)$ is equal to the connecting graph of
$\{(k,c_k)\}_{0\leq k\leq d}$ on $\R^2$.

\noindent \fbox{Part 1. Newton polygon of $\alpha_1$ over
$\Omega'_a$.} {}From now on let $\bM$ be the (infinite) matrix
representing the $\alpha_1$ action on $\Omega'_a$-space
$\cH(\bA_r)$ with respect to the basis $\vec{b}_{\rm w}$. (See
Table 1.) Write
\begin{equation}\label{E:detM}
\det(1-T\bM)=1+\sum_{k=1}^{\infty}C_kT^k\in\cO'_a[[T]].
\end{equation}
Take the minimal $p$-adic valuation of all entries in each row,
and put them in a nondecreasing order. For any $i\geq 1$ let
$m_i(\bM)$ denote the $i$-th smallest row $p$-adic valuation of
$\bM$ (counting multiplicity). For every $k\geq 1$ let
$c_k(\bM):=\sum_{i=1}^{k}m_i(\bM)$. By Theorem \ref{T:reduce1},
one has
\begin{equation}\label{E:C}
\ord_pC_{J_1,J}^{n,i}\geq
\left\{
\begin{array}{ll}
 \frac{n}{d_{J_1}}   & \mbox{for $J_1=1,2$ and $n,i\geq 0$}\\
 \frac{n-1}{d_{J_1}} & \mbox{for $J_1\geq 3$ and $n\geq 1, i\geq 0$}.
\end{array}
\right.
\end{equation}
This implies that $m_i(\bM)\geq m_i$ for all $1\leq i\leq
d-\ell+1$. Thus by arithmetic of Newton polygons (see \cite[Lemma
1.6]{Dwork1}) and Fredholm theory (see \cite[Proposition 7 and its
proof] {Serre}) one has that $\NP_q(\det(1-T\bM)\bmod
T^{d-\ell+2})$ lies above the connecting graph of $\{(k,c_k)\in
\R^2\}_{0\leq k\leq d-\ell+1}$. The latter is precisely $\HP(\A)$
as remarked earlier.

\noindent\fbox{Part 2. Newton polygon of $\alpha_1$ over
$\Omega'_1$.} By the normal basis theorem, there exists
$\xi\in\Omega_a$ such that $\vec{\xi}:=\{\xi^{\tau^t}\}_{0\leq
t\leq a-1}$ is a basis for $\Omega'_a$ over $\Omega'_1$. Let $\bN$
be the (infinite) matrix representing $\alpha_1$ with respect to
the basis $\vec{b}_{{\rm w},\Omega'_1}$ for $\cH(\bA_r)$ as
$\Omega'_1$-space, where $\vec{b}_{{\rm w},\Omega'_1}$ consists of
$Z_j^i \xi^{\tau^t}$ for $1\leq j\leq \ell$, $0\leq t\leq a-1$ and
$i\geq 0$ (where $i=0$ only if $j=1$). Write $\det(1-T\bN)
=1+\sum_{k=1}^{\infty}D_kT^k\in\cO'_1[[T]].$ We have $
\alpha_1(Z_{J}^i\xi^{\tau^t}) =\alpha_1(Z_{J}^i)\xi^{\tau^{t-1}} =
\sum_{J_1=1}^{\ell}\sum_{n=0}^{\infty}
   C_{J_1,J}^{n,i} Z_{J_1}^n \xi^{\tau^{t-1}}.
$
Recall the lower bound of $\ord_p(C_{J_1,J}^{n,i})$ given in
(\ref{E:C}). In the two sequences
\begin{eqnarray}
\label{E:seq1}
&&\{m_1(\bN),m_2(\bN),\cdots,\cdots, m_{a(d-\ell+1)}(\bN)\} \\
\label{E:seq2}
&&\{\underbrace{m_1,\cdots, m_1}_{a},\cdots,
\underbrace{m_{d-\ell+1},\cdots, m_{d-\ell+1}}_{a}\},
\end{eqnarray}
one notes that (\ref{E:seq1}) {\em dominates} (\ref{E:seq2}) in
the sense that the $i$-th term of the former sequence is greater
than or equal to that of the latter. Thus $1/a\cdot
\NP_p(\det_{\Omega_1}(1-T\alpha_1) \bmod T^{a(d-\ell+1)+1})$ lies
above the connecting graph of $\{(k,c_k)\}_{0\leq k \leq
d-\ell+1}$, that is $\HP(\A)$. By applying  Proposition
\ref{P:reduce1}, one now concludes that $\NP(\bar{f};\F_q)$ lies
over $\HP(\A)$.

\noindent\fbox{Part 3. Newton and Hodge coincide if and only if
$p\equiv 1\bmod (\lcm\; d_j)$.} One notes that after permuting our
basis $\vec{b}_{\rm w}$ for $\cH(\bA_r)$ we can arrive at a matrix
$\bM$ of $\alpha_1$ in block form satisfying the hypothesis of
\cite[Theorem 2.4]{Wan:1}. Let $\bM_a$ be the matrix representing
$\alpha_a$ over $\Omega'_a$, then one knows that
$\bM_a=\bM\bM^{\tau^{-1}}\cdots \bM^{\tau^{-(a-1)}}$. Note that
$\NP_q(\det(1-T\bM_a))$ and $\HP(\A)$ meet at the point with
$x$-coordinate $d-\ell+1$ (see Proposition \ref{P:reduce1}). Let
$\HP(\A)_{<1}$ denote the slope $<1$ part of $\HP(\A)$, which has
horizontal length $d-\ell+1$. By \cite[Theorem 2.4, Corollary
2.5]{Wan:1}, one can show that $\NP_q(\det(1-T\bM_a)\bmod
T^{d-\ell+2}) =\HP(\A)_{<1}$ if and only if
$\NP_p(\det(1-T\bM)\bmod T^{d-\ell+2})=\HP(\A)_{<1}$. That is, it
is enough to show $\NP_p(\det(1-T\bM)\bmod
T^{d-\ell+2})=\HP(\A)_{<1}$.

Let $\bM_{<1}$ be the principle submatrix of $\bM$ consisting of
all $C_{J_1,J}^{n,i}$ with
$$
\left\{
\begin{array}{ll}
 0\leq n\leq d_1-1        &\mbox{for $J_1=1$};\\
 1\leq n\leq d_2-1        &\mbox{for $J_1=2$};\\
 1\leq n\leq d_{J_1}  &\mbox{for $J_1\geq 3$}.
\end{array}
\right.
$$
One notices that $\bM_{< 1}$ has $d-\ell+1$ rows in total. By
(\ref{E:C}), every row of $\bM$ outside these $d-\ell+1$ rows has
its minimal $p$-adic valuation greater than or equal to $1$.
{}From matrix arithmetic
of Fredholm theory, it is not hard to conclude that all segments
of $\NP_p(\det(1-T\bM))_{<1}$ have to ``come from"
$\det(1-T\bM_{<1})$ in the following sense. In (\ref{E:detM}) let
$t$ be the biggest integer such that $\NP_q(\sum_{k=0}^{t}C_kT^k)$
has all slopes less than $1$,
then for all $k\leq n$ one has $C_k=\sum_N \pm
\det N$ where $N$ ranges over all $k\times k$ principal submatrices
of $\bM_{<1}$.

{}Now we assume that $p\equiv 1\bmod (\lcm\; d_j)$ where $j$
ranges from $1$ to $\ell$. By Remark \ref{R:main2}, the slope-$0$
segment of the Hodge polygon is always achieved. This saves us
{}from considering the corresponding rows in $\bM_{<1}$. By
Theorem \ref{T:reduce1}, for $J_1=1,2$ (resp. $J_1\geq 3$) one has
that $C_{J_1,J}^{n,i}$ in the submatrix $\bM_{<1}$ achieves its
minimal row $p$-adic value $\frac{n}{d_{J_1}}$ (resp.,
$\frac{n-1}{d_{J_1}}$) uniquely at $n=i$. In summary, all minimal
row $p$-adic valuations are achieved uniquely on the diagonal of
$\bM_{<1}$. These minimal row $p$-adic valuations are precisely
the rational numbers less than $1$ listed in (\ref{E:HP}). By
arithmetic of Fredholm theorem and analysis in Part 1,
$\NP_p(\det(1-T\bM)\bmod T^{d-\ell+2}) =\NP_p(\det(1-T\bM_{<1}))
=\HP(\A)_{<1}$. Combining with the above paragraph, we have shown
that $\NP_q(\det(1-T\bM_a)\bmod T^{d-\ell+2}) = \HP(\A)_{<1}$. By
Proposition \ref{P:reduce1}, one concludes that
$\NP(\bar{f};\F_q)=\HP(\A)$.

Conversely, suppose $\NP(\bar{f};\F_q)$ coincides with Hodge. By
the above argument, one has $m_i(\bM)=m_i$ for all $i$. Because
every minimal row $p$-adic valuation is achieved only at $J=J_1$
(it lies in the diagonal blocks in the matrix $\bM$), the Newton
polygon of $\bM$ lies above the end-to-end join of those of
$\bM_j$ for $1\leq j\leq \ell$ where
$\bM_j:=\{C_{j,j}^{n,i}\}_{1\leq n,i\leq d_j-1}$. Thus the Newton
polygon of $\bM_j$ has to coincide with its Hodge (its Hodge
polygon is defined in the obvious sense). By the remark in
the second last paragraph in Section 1,
one can shift the pole to $\infty$ so that we may assume
$p\not\equiv 1\bmod d_1$. Since $\ord_p(C_{1,1}^{n,i})=n/d_1$ for
some $1\leq i_n\leq d_1-1$ for all $1\leq n\leq d_1-1$, by Theorem
\ref{T:sharp}, the condition $\cC(np-i_n,d_1)$ holds, and it is
$np\equiv i_n\bmod d_1$. Since $p\not\equiv 1\bmod d_1$, one has
$n\neq i_n$ for every $n$. From simple linear algebra, one sees
that the first slope of the Newton polygon of $\bM_j$ is
greater than or equal to $\ord_p(C_{1,1}^{1,i_1})>1/d_1$. A contradiction.

This finishes the proof of Theorem \ref{T:main1}.

\begin{acknowledgments}
I am grateful to Bjorn Poonen and Steve Sperber for communicating
to me of their conjecture and for stimulating discussions. I thank
Hanfeng Li for many exchanges of ideas and in particular for
constructive comments on Section \ref{S:2.3}.
\end{acknowledgments}


\begin{thebibliography}{99}

\bibitem{Berthelot}
{\sc Pierre Berthelot:} Cohomologie rigide et th\'eorie de Dwork:
le cas des sommes exponentielles[Rigid cohomology and
Dwork theory: the case of exponential sums],
in {\it Cohomologie p-adique}
Soci\'et\'e Math\'ematique de France,
Ast\'erisque {\bf 119--120}(1984), 17-49 (French).

\bibitem{Berthelot:2}
{\sc Pierre Berthelot:}
 G\'eom\'etrie rigide et cohomologie des vari\'et\'es alg\'ebriques
 de caract\'eristique $p$ [Rigid geometry and cohomology of algebraic
 varieties in characteristic $p$].
 Introductions aux cohomologies $p$-adiques
 (Luminy, 1984). M\'em. Soc. Math. Fr.(N.S.) {\bf 23} (1986), 7--32 (French).

\bibitem{Bombieri}
{\sc Enrico Bombieri:}
On exponential sums in finite fields,
{\it Amer. J. Math.} {\bf 88} (1966), 71--105.

\bibitem{BGR}
{\sc S. Bosch, U. Guntzer, R. Remment:}
Non-Archimedean analysis,
{\it Grundlehren der mathematischen Wissenschaften}
Vol. {\bf 261}, Springer-Verlag, Berlin, 1984.

\bibitem{Crew}
{\sc Richard Crew:}
Etale $p$-covers in characteristic $p$,
Compositio Math., {\bf 52} (1984), no.1,
31--45.

\bibitem{Dwork1}
{\sc Bernard Dwork:}
On the zeta function of a hypersurface.
{\it Inst. Hautes \'Etudes Sci. Publ. Math.}{\bf 12}
(1962), 5--68.

\bibitem{Dwork2}
{\sc Bernard Dwork:}
On the zeta function of a hypersurface. II.
{\it Ann. of Math.(2)}
{\bf 80} (1964), 227--299.

\bibitem{Dwork3}
{\sc Bernard Dwork:}
{\it Lectures on $p$-adic differential equations.}
{Grundlehren der mathematischen Wissenschaften} {\bf 253},
Springer-Verlag, New York, 1982.

\bibitem{Katz}
{\sc N. Katz:} Slope filtration of F-crystals, {\it Ast\'erisque}
{\bf 63} (1979), 113--163.

\bibitem{Koblitz:0}
{\sc Neal Koblitz:}
$p$-adic variation of the zeta function over families of
varieties defined over finite fields,
{\em Comp. Math.}
{\bf 31}
(1975), no.2,
119--218.

\bibitem{Koblitz}
{\sc Neal Koblitz:}
$p$-adic numbers, p-adic analysis, and zeta-functions,
(Second edition),
{\it Graduate Texts in Mathematics} {\bf 58}.
Springer-Verlag, New York, 1984.

\bibitem{Krasner}
{\sc Marc Krasner:}
Prolongement analytique uniforme et multiforme dans les corps
valu\'es complets,
{\it Les Tendances G\'eom.
en Alg\`ebre et Th\'eorie des Nombres}, \'Editions du Centre National de la
Recherche Scientifique, CNRS, Paris, 1966, pp. 97--141 (French).

\bibitem{Li-Zhu:1}
{\sc Hanfeng Li and Hui June Zhu:} Asymptotic variation of sheaves
of one-variable exponential sums. 2003, preprint,
{\tt http://arXiv.org/abs/math.NT/0312423}.

\bibitem{Mazur}
{\sc Barry Mazur:} Frobenius and the Hodge filtration, Bull.
A.M.S. {\bf 78} (1972), 653--667.

\bibitem{Monsky-Washnitzer:I}
{\sc P. Monsky; Washnitzer:}
Formal cohomology I,
{\it Ann. of Math.}
{\bf 88} (1968), 181--217.

\bibitem{Monsky:II}
{\sc P. Monsky:} Formal cohomology II, {\it Ann. of Math.}
{\bf 88} (1968), 218--238.

\bibitem{Monsky:III}
{\sc P. Monsky:}
Formal cohomology III, {\it Ann. of Math.}
{\bf 93} (1971), 315--343.

\bibitem{Reich}
{\sc Daniel Reich:}
A $p$-adic fixed point formula.
{\it Amer. J. Math.}
{\bf 91} (1969), 835--850.

\bibitem{Robba:0}
{\sc Philippe Robba:}
Fonctions analytiques sur les corps valu\'es untram\'etriques
complets, in
{\it Prolongement analytique et alg\`ebres banach ultram\,etriques},
Ast\'erisque {\bf 10},
Soci\'et\'e Math\'ematique de France, Paris, 1973, pp. 109--220 (French).

\bibitem{Robba:1}
{\sc Philippe Robba:}
Index of p-adic differential operators III. Application to twisted
exponential sums, in
{\it Cohomologie p-adique}
Ast\'erisque {\bf 119--120}
(1984), 191--266.

\bibitem{Robba:2}
{\sc Philippe Robba:}
Une Introduction n\"aive aux cohomologies de Dwork, {\em Soc.
Math. Fr. M\'emoire (N.S.)} {\bf 23} (1986), 61--105 (French).

\bibitem{Schikhop}
{\sc W.H. Schikhof:} Ultrametric calculus. an introduction to
$p$-adic analysis. {\it Cambridge Studies in Advanced
mathematics}, {vol.\bf 4}. Cambridge University Press, Cambridge,
1984.

\bibitem{Serre}
{\sc J-P Serre:}
Endomorphismes compl\`etements continus des espaces de Banach
$p$-adique,
{\it Inst. Hautes. \'Etudes Sci. Publ. Math.} {\bf 12}
(1962), 69--85 (French).

\bibitem{Sperber:0}
{\sc Steven Sperber:}
Congruence properties of the hyper-Kloosterman sum.
{\it Compositio Math.} {\bf 40} (1980), no.1,
3--33.

\bibitem{Sperber}
{\sc Steven Sperber:} On the $p$-adic theory of exponential sums,
{\it Amer. J. Math.,} {\bf 109} (1986), no.2,
255--296.

\bibitem{Geer-Vlugt:92}
{\sc Gerard van der Geer; Marcel van der Vlugt:}
Reed-Muller codes and supersingular curves. I.
{\it Compositio Math.} {\bf 84}
(1992), no. 3,
333--367.

\bibitem{Geer-Vlugt:95}
{\sc Gerard van der Geer; Marcel van der Vlugt:}
On the existence of supersingular curves of given genus.
{\it J. Reine Angew. Math.}
{\bf 458} (1995), 53--61.

\bibitem{Wan:1}
{\sc Daqing Wan:} Newton polygons of zeta functions and L
functions. {\it Ann. of Math.}, {\bf 137} (1993), 249--293.

\bibitem{Zhu:1}
{\sc Hui June Zhu:}
$p$-adic variation of $L$ functions of one variable exponential
sums, I. {\it Amer. J. Math.} {\bf 125} (2003).

\bibitem{Zhu:2}
{\sc Hui June Zhu:}
Asymptotic variation of $L$ functions of one-variable exponential
sums.
{\it J. Reine Angew. Math.} {\bf 572} (2004),
219--233.


\end{thebibliography}
\end{document}